\documentclass{amsart}
\usepackage{amssymb}
\usepackage{latexsym}
\newcommand{\Z}{{\mathbb Z}}
\newcommand{\R}{{\mathbb R}}
\newcommand{\C}{{\mathbb C}}
\newcommand{\D}{{\mathbb D}}

\newtheorem{theorem}{Theorem}

\newtheorem{lemma}{Lemma}[section]
\newtheorem{prop}[lemma]{Proposition}

\let\Re=\undefined\DeclareMathOperator*{\Re}{Re}

\newcounter{smalllist}

\makeatletter\renewcommand{\thesmalllist}{\@alph\c@smalllist}\makeatother

\begin{document}

\title{Verblunsky Coefficients With Coulomb-Type Decay}

\author{David Damanik}

\address{Mathematics 253-37, California Institute of Technology, Pasadena, CA 91125, USA}

\email{damanik@caltech.edu}

\date{\today}

\begin{abstract}
We show that probability measures on the unit circle associated with Verblunsky
coefficients obeying a Coulomb-type decay estimate have no singular continuous component.
\end{abstract}

\maketitle

\begin{center}
\textit{Dedicated to Barry Simon on the occasion of his 60th birthday.}
\end{center}

\bigskip

%%%%%%%%%%%%%%%%%%%%%%%%%%%%%%%%%%%%%%%%%%%%%%%%%%%%%%%%%%%%%%%%%%%%%%%%%%%%%%%%%%%%%
\section{Introduction}
%%%%%%%%%%%%%%%%%%%%%%%%%%%%%%%%%%%%%%%%%%%%%%%%%%%%%%%%%%%%%%%%%%%%%%%%%%%%%%%%%%%%%

Let $d\mu$ be a probability measure on $\R / (2 \pi \Z)$ that is not supported on a
finite number of points. Then, using the Gram-Schmidt procedure, we may find polynomials
$\varphi_n(z)$ that obey
$$
\int_0^{2\pi} \overline{\varphi_m(e^{i\eta})} \varphi_n(e^{i\eta}) \, d\mu(\eta) =
\delta_{m,n}.
$$
We also consider the monic orthogonal polynomials $\Phi_n(z)$. They obey the Szeg\H{o}
recursion
$$
\Phi_{n+1}(z) = z \Phi_n(z) - \overline{\alpha}_{n} \Phi_n^*(z),
$$
where $\Phi_n^*(z) = z^n \overline{\Phi_n(1/\overline{z})}$. The $\alpha_n$ are called
Verblunsky coefficients and they belong to the unit disk $\D = \{ z \in \C : |z| < 1 \}$.
Conversely, every $\alpha \in \times_{n=0}^\infty \D$ corresponds to a unique measure.
See \cite{simon1,simon,szego} for background material on orthogonal polynomials on the
unit circle (OPUC).

In this paper we are interested in the measures associated with Verblunsky coefficients
that have Coulomb-type decay. To motivate our study, let us recall the following result
of Golinskii and Ibragimov \cite{gi}:
$$
\sum_{n = 0}^\infty (n+1) |\alpha_n|^2 < \infty \quad \Rightarrow \quad d\mu_{{\rm sing}}
= 0.
$$
Here, $d\mu_{{\rm sing}} = d\mu_{{\rm sc}} + d\mu_{{\rm pp}}$, where $d\mu = d\mu_{{\rm
ac}} + d\mu_{{\rm sc}} + d\mu_{{\rm pp}}$ is the Lebesgue decomposition of $d\mu$ into an
absolutely continuous (with respect to Lebesgue measure) piece, a singular continuous
piece, and a pure point piece.

The natural class of Verblunsky coefficients having true Coulomb decay, that is,
$|\alpha_n| = O(1/n)$, is outside the scope of the result above. More generally, one may
be interested in the class of Verblunsky coefficients satisfying
\begin{equation}\label{vcass}
\sum_{n = 0}^N (n+1) |\alpha_n|^2 \le A \log N
\end{equation}
for some $A < \infty$. The following extension, due to Simon \cite{s}, of the result of
Golinskii and Ibragimov covers a portion of this class:
\begin{equation}\label{simres}
\alpha \text{ satisfies \eqref{vcass} for some } A < \tfrac14 \quad \Rightarrow \quad
d\mu_{{\rm sing}} = 0.
\end{equation}
Simon also shows that for every $A > \frac{1}{4}$, there is an example satisfying
\eqref{vcass} with $d\mu_{{\rm pp}} \not=0$. Thus, the result \eqref{simres} is almost
sharp. The latter result is an OPUC analogue of the classical Wigner-von Neumann example
that exhibits an embedded eigenvalue for a half-line Schr\"odinger operator with $O(1/x)$
potential \cite{WvN}.

The pure point component is further studied in \cite{simon}. There it is shown (see
Theorem~10.12.7) that if \eqref{vcass} holds for some $A$, then $d\mu$ has at most $K$
pure points, where $K$ is the unique integer with $K \le 4A < K + 1$. Following this
theorem, Simon writes that it is an intriguing open question if \eqref{vcass} implies
$d\mu_{{\rm sc}} = 0$. There are two reasons why one expects a positive answer to this
question. Intuitively, it should be easier to have infinitely many pure points than a
singular continuous component, so that the result just quoted supports the conjecture
that a singular continuous piece should be impossible. On the other hand, Kiselev has
proven the absence of singular continuous spectrum for half-line Schr\"odinger operators
with $O(1/x)$ potentials \cite{k}.

Our goal here is to give an affirmative answer to Simon's question and prove the
following theorem:

\begin{theorem}\label{main}
Suppose there is $A < \infty$ such that $\alpha$ satisfies \eqref{vcass}. Then,
$d\mu_{{\rm sc}} = 0$.
\end{theorem}

Since it is also shown in \cite{s} that \eqref{vcass} with $A = 1/4$ implies $d\mu_{{\rm
pp}} = 0$, it follows from Theorem~\ref{main} that \eqref{simres} may be strengthened to
$$
\alpha \text{ satisfies \eqref{vcass} for some } A \le \tfrac14 \quad \Rightarrow \quad
d\mu_{{\rm sing}} = 0,
$$
which is optimal by the discussion above.

The overall strategy in our proof of Theorem~\ref{main} will be inspired by Kiselev
\cite{k}. This will require some preparatory work. We first recall Pr\"ufer variables and
the Bernstein-Szeg\H{o} Approximation to $d\mu$ in Section~2 and prove a comparison lemma
which is related to the Chebyshev-Markov Moment Problem. Then, we consider the support of
$d\mu_{{\rm sing}}$ in Section~3 and prove that it has Hausdorff dimension zero. This is
a result in the spirit of Remling \cite{r2} who proved results of this flavor for
half-line Schr\"odinger operators. Finally, we prove Theorem~\ref{main} in Section~4 by
working out the OPUC analogue of Kiselev's ideas from \cite{k}.

\bigskip

\noindent\textit{Acknowledgments.} The author would like to thank Rowan Killip and
Christian Remling for useful conversations.

%%%%%%%%%%%%%%%%%%%%%%%%%%%%%%%%%%%%%%%%%%%%%%%%%%%%%%%%%%%%%%%%%%%%%%%%%%%%%%%%%%%%%
\section{Pr\"ufer Variables and Bernstein-Szeg\H{o} Approximation}
%%%%%%%%%%%%%%%%%%%%%%%%%%%%%%%%%%%%%%%%%%%%%%%%%%%%%%%%%%%%%%%%%%%%%%%%%%%%%%%%%%%%%

Let $\{\alpha_n\}$ be the Verblunsky coefficients of a nontrivial probability measure
$d\mu$ on $\partial \D$. As mentioned above, the $\alpha$'s give rise to a sequence
$\{\Phi_n(z)\}$ of monic polynomials (via the Szeg\H{o} recursion) that are orthogonal
with respect to $d\mu$. For $\beta \in [0,2\pi)$, we also consider the monic polynomials
$\{\Phi_n(z,\beta)\}$ that are associated in the same way with the Verblunsky
coefficients $\{ e^{i\beta} \alpha_n \}$.

Let $\eta \in [0,2\pi)$. Define the Pr\"ufer variables by
$$
\Phi_n(e^{i\eta},\beta) = R_n(\eta,\beta) \exp \left[ i(n \eta + \theta_n(\eta,\beta))
\right],
$$
where $R_n > 0$, $\theta_n \in [0,2\pi)$, and $|\theta_{n+1} - \theta_n| < \pi$. These
variables obey the following pair of equations:
\begin{align*}
\frac{R_{n+1}^2(\eta,\beta)}{R_n^2(\eta,\beta)} & = 1 + | \alpha_n |^2 - 2 \Re
\left(\alpha_n e^{i[(n+1)\eta + \beta + 2
\theta_n(\eta,\beta)]} \right), \\
e^{-i(\theta_{n+1}(\eta,\beta) - \theta_n(\eta,\beta))} & = \frac{1 - \alpha_n
e^{i[(n+1)\eta + \beta + 2 \theta_n(\eta,\beta)]}}{\left[1 + | \alpha_n |^2 - 2 \Re
\left(\alpha_n e^{i[(n+1)\eta + \beta + 2 \theta_n(\eta,\beta)]} \right)\right]^{1/2}}.
\end{align*}
We also define $r_n(\eta,\beta) = |\varphi_n(\eta,\beta)|$.

When $\{ \alpha_n \} \in \ell^2$,
\begin{equation}\label{fs}
r_n(\eta,\beta) \sim R_n(\eta,\beta) \sim \exp \left( - \sum_{j=0}^{n-1} \Re (\alpha_j
e^{i[(j+1)\eta + \beta + 2 \theta_j(\eta,\beta)]}) \right).
\end{equation}
(We write $f_n \sim g_n$ if there is $C > 1$ such that $C^{-1}g_n \le f_n \le C g_n$ for
all $n$.) For the Pr\"ufer equations and \eqref{fs}, see \cite[Theorems~10.12.1 and
10.12.3]{simon}.

Next we recall the Bernstein-Szeg\H{o} Approximation of $d\mu$. The measure $d\mu_n$
associated with Verblunksy coefficients $\alpha_0,\ldots,\alpha_{n-2},
\alpha_{n-1},0,0,\ldots$ is given by
\begin{equation}\label{bsa}
d\mu_n(\eta) = \frac{d\eta}{2\pi  r_n^2(\eta,0)};
\end{equation}
compare \cite[Theorem~1.7.8]{simon1}.

If $d\mu$ and $d\nu$ are two measures whose first $n$ Verblunsky coefficients coincide
(i.e., $\alpha_k(d\mu) = \alpha_k(d\nu)$, $0 \le k \le n-1$), their moments up to order
$n$ are the same (see, e.g., \cite[Theorem~1.5.5.(ii)]{simon1}). Consequently, given a
Laurent polynomial, $f(\eta) = \sum_{k = -n}^n f_k e^{ik\eta}$, we have
\begin{equation}\label{bsaint}
\int_0^{2\pi} f(e^{i\eta}) \, d\mu(\eta) = \int_0^{2\pi} f(e^{i\eta}) \, d\nu(\eta).
\end{equation}

\begin{lemma}\label{bsaest}
Suppose $d\mu$ and $d\nu$ are two measures whose first $n$ Verblunsky coefficients
coincide. For every $\kappa > 0$, $n \in \Z_+$, and every interval $I \subseteq \partial
\D$ of length $\delta \ge n^{-1/(2 + \kappa)}$, we have
\begin{equation}\label{munu}
\mu(I) \le \nu(3I) + C \delta^\kappa.
\end{equation}
\end{lemma}

\noindent\textit{Remarks.} (a) In \eqref{munu}, $3I$ denotes the interval of length
$3\delta$ that has the same center as $I$ and $C$ is a constant that depends only on
$\kappa$. Alternatively, one may choose a universal $C$ for which \eqref{munu} holds for
all $\delta \le \delta_0$ (and hence $n \ge n_0(\kappa)$).\\[1mm]
(b) Since the estimate \eqref{munu} is sufficient for our purpose and has a short and
elementary proof, we content ourselves with this explicit statement. We do want to point
out, however, that it is closely related to the Chebyshev-Markov Moment Problem: If we
fix $n$ initial moments and an interval $I$, what are the extremal values of $\mu(I)$
when $\mu$ ranges over all measures that have the prescribed moments? A wealth of
material dealing with this problem may be found, for example, in \cite{krein,kreinnudel}.
An analogue of these classical results for Schr\"odinger operators in $L^2(0,\infty)$ was
recently found in \cite{rem2}.

\begin{proof}
Without loss of generality, we assume that $I = (-\frac{\delta}{2},\frac{\delta}{2})$.
Consider the Fej\'er kernel,
$$
F_n(\eta) = \sum_{k=-n}^n \left( 1 - \frac{|k|}{m+1} \right) e^{ik\eta} = \frac{1}{n+1}
\left( \frac{\sin \tfrac{n+1}{2} \eta}{\sin \tfrac{1}{2} \eta} \right),
$$
and let
$$
\sigma_n(\eta) = (F_n * \chi_{2I})(\eta) = \frac{1}{2\pi} \int_0^{2\pi} F_n(\tau)
\chi_{2I}(\eta - \tau) \, d\tau.
$$
Clearly,
\begin{equation}\label{fejer1}
|\sigma_n (\eta)| \le 1 \text{ for all } \eta.
\end{equation}
Moreover, by \eqref{bsaint}, it follows that
\begin{equation}\label{fejer2}
\int_0^{2\pi} \sigma_n(\eta) \, d\mu(\eta) = \int_0^{2\pi} \sigma_n(\eta) \, d\nu(\eta).
\end{equation}
Note that
$$
\sigma_n(\eta) - \chi_{2I}(\eta) = \frac{1}{\pi} \int_0^{\pi} F_n(\tau) \left[
\frac{\chi_{2I}(\eta - \tau) + \chi_{2I}(\eta + \tau)}{2} - \chi_{2I}(\eta) \right] \,
d\tau.
$$
When $||\eta| - \delta | \ge \frac{\delta}{2}$, this gives
$$
\sigma_n(\eta) - \chi_{2I}(\eta) = \frac{1}{\pi} \int_{\frac{\delta}{2}}^\pi F_n(\tau)
\left[ \frac{\chi_{2I}(\eta - \tau) + \chi_{2I}(\eta + \tau)}{2} - \chi_{2I}(\eta)
\right] \, d\tau.
$$
Consequently, for these values of $\eta$, we have
$$
\left| \sigma_n(\eta) - \chi_{2I}(\eta) \right| \le \frac{2}{\pi}
\int_{\frac{\delta}{2}}^\pi F_n(\tau) \, d\tau \le \frac{2}{n+1} \frac{1}{\sin^2
\frac{\delta}{2}} \lesssim \frac{1}{\delta^2 n} \le \delta^\kappa,
$$
where we used the assumption $\delta \ge n^{-1/(2 + \kappa)}$ in the last step. Thus,
\begin{equation}\label{fejer3}
|\sigma_n(\eta) - \chi_{2I}(\eta)| \lesssim \delta^\kappa \text{ for all $\eta$
satisfying } ||\eta| - \delta | \ge \frac{\delta}{2}.
\end{equation}
The assertion of the lemma is an immediate consequence of \eqref{fejer1}--\eqref{fejer3}.
\end{proof}

%%%%%%%%%%%%%%%%%%%%%%%%%%%%%%%%%%%%%%%%%%%%%%%%%%%%%%%%%%%%%%%%%%%%%%%%%%%%%%%%%%%%%
\section{Zero-Dimensionality of the Singular Part}
%%%%%%%%%%%%%%%%%%%%%%%%%%%%%%%%%%%%%%%%%%%%%%%%%%%%%%%%%%%%%%%%%%%%%%%%%%%%%%%%%%%%%

In this section we show that the singular part of $d\mu$ must be supported on a set of
zero Hausdorff dimension if the Verblunsky coefficients obey \eqref{vcass}. Results of
this kind were obtained in the context of Schr\"odinger operators by Remling \cite{r2},
Christ-Kiselev \cite{ck}, and Damanik-Killip \cite{dk}, for example. We will follow ideas
from \cite{dk} rather closely.

\begin{prop}\label{p31}
Assume \eqref{vcass}. Then the set
$$
S = \{ \eta \in [0,2\pi) : R_n(\eta,\beta) \text{ is unbounded for some } \beta \}
$$
has zero Hausdorff dimension. Consequently, $d\mu_{{\rm sing}}$ is supported on a set of
zero Hausdorff dimension.
\end{prop}

Clearly, \eqref{vcass} implies $\{\alpha_n\} \in \ell^2$. Therefore, because of
\eqref{fs}, our goal is to show that
$$
A(n,\eta,\beta) = \sum_{j=0}^{n-1} \alpha_j e^{i[(j+1)\eta + \beta + 2
\theta_j(\eta,\beta)]}
$$
is a bounded function of $n$ for all $\beta$, provided that $\eta$ is away from a set of
zero Hausdorff dimension.

\begin{lemma}\label{klstool}
If
$$
\hat{\alpha}(\eta,n) = \lim_{N \to \infty} \sum_{j = n}^N \alpha_j e^{i j \eta}
$$
exists and obeys
\begin{equation}\label{fsr}
\sum_{j=1}^\infty | \hat{\alpha} (\eta,j) \alpha_{j-1} | < \infty,
\end{equation}
then $\eta \not\in S$.
\end{lemma}

\begin{proof}
We will show that $A(n,\eta,\beta)$ is bounded (in $n$) for every $\beta \in [0,2\pi)$
when \eqref{fsr} holds. The assertion then follows from \eqref{fs}.

Write $\gamma_j(\eta,\beta) = (j+1)\eta + \beta + 2 \theta_j(\eta,\beta)$. We have
\begin{align*}
A(n,\eta,\beta) & = \sum_{j=0}^{n-1} \left[\hat{\alpha} (\eta,j) -
\hat{\alpha} (\eta,j+1)\right] e^{i\gamma_j(\eta,\beta) - i j \eta} \\
& = \sum_{j=1}^{n-1} \hat{\alpha} (\eta,j) \left[ e^{i\gamma_j(\eta,\beta)} -
e^{i(\gamma_{j-1}(\eta,\beta) + \eta)} \right] e^{-i j \eta} + O(1).
\end{align*}
Since
\begin{align*}
| e^{i\gamma_j(\eta,\beta)} - e^{i(\gamma_{j-1}(\eta,\beta) + \eta)}| & \le |
\gamma_j(\eta,\beta) - \gamma_{j-1}(\eta,\beta) - \eta | \\ & = 2 | \theta_j(\eta,\beta)
- \theta_{j-1}(\eta,\beta) | \\ & \lesssim |\alpha_{j-1}|,
\end{align*}
boundedness of $A(n,\eta,\beta)$ follows.
\end{proof}

\begin{lemma}\label{SZ}
Let $d\nu$ be a positive measure on $[0,2\pi)$. For each $\varepsilon\in(0,1)$ and every
measurable function $m:[0,2\pi) \to \Z_+$,
$$
\Biggl\{ \int\ \Biggl| \sum_{n=0}^{m(\eta)}  c_n e^{-in\eta} \Biggr| \, d\nu(\eta)
\Biggr\}^2 \lesssim \mathcal{E}_\varepsilon (\nu) \sum_{n=0}^\infty (n+1)^{1-\varepsilon}
\big|c_n\big|^2\! ,
$$
where $\mathcal{E}_\varepsilon$ denotes the $\varepsilon$-energy of $d\nu$:
$\mathcal{E}_\varepsilon (\nu) = \int \int (1+|x-y|^{-\varepsilon}) \, d\nu(x) \,
d\nu(y)$.
\end{lemma}

\begin{proof}
This follows by slightly adjusting the calculation from \cite[\S XIII.11,
p.~196]{Zygmund} (see also \cite[\S V.5]{Carleson}).
\end{proof}

\begin{proof}[Proof of Proposition~\ref{p31}.]
We will apply the criterion of Lemma~\ref{klstool}. Let us first note that by the theorem
of Salem-Zygmund, the series defining $\hat{\alpha}$ converges off a set of zero
Hausdorff dimension. Therefore, we may exclude from consideration those values of $\eta$
for which $\hat{\alpha}$ is not defined.

By applying the Cauchy-Schwarz inequality to dyadic blocks, for example, we see that
\eqref{vcass} implies $n^{-\varepsilon/4} \alpha_n \in \ell^1$ for all $\varepsilon > 0$.
Hence the proposition will follow from Lemma~\ref{klstool} once we prove that for all
$\varepsilon > 0$, the set of $\eta$ for which $n^{\varepsilon/4} \hat{\alpha}(\eta,n)$
is unbounded is of Hausdorff dimension no more than $\varepsilon$.

Let $m(\eta)$ be a measurable integer-valued function on $[0,2\pi)$. Because of
\eqref{vcass}, Lemma~\ref{SZ} implies
\begin{align*}
   \int\ \Biggl| \sum_{n=m_l(\eta)}^{2^{l+1} - 1} \alpha_n e^{i n \eta} \Biggr| \, d \nu(\eta)
&= \int\ \Biggl| \sum_{n=0}^{\tilde{m}_l(\eta)} \alpha_{2^{l+1}-1-n} e^{-i n \eta} \Biggr| \, d \nu(\eta) \\
&\lesssim \Biggl\{ \sum_{n=2^l}^{2^{l+1}-1} (n+1)^{1-\varepsilon} \big|\alpha_n\big|^2
\Biggr\}^{1/2} \sqrt{\mathcal{E}_\varepsilon (\nu)} \\
&\lesssim \sqrt{l} \, 2^{-\varepsilon l/2} \sqrt{\mathcal{E}_\varepsilon (\nu)}
\end{align*}
where $m_l(\eta)=\max\{m(\eta),2^l\}$, $\tilde{m}_l(\eta) = \min \{2^l-1, 2^{l+1} - 1 -
m(\eta)\}$, and sums with lower index greater than their upper index are to be treated as
zero. Multiplying both sides by $2^{\varepsilon l/4}$, summing this over $l$, and
applying the triangle inequality on the left gives
$$
\int\ \Biggl| m(\eta)^{\varepsilon/4} \sum_{n=m(\eta)}^{\infty} \alpha_n e^{i n \eta}
\Biggr| \, d \nu(\eta) \lesssim \sqrt{\mathcal{E}_\varepsilon (\nu)}.
$$
That is, for any measurable integer-valued function $m(\eta)$,
$$
\int m(\eta)^{\varepsilon/4} \bigl| \hat{\alpha}(\eta, m(\eta))\bigr| \, d\nu \lesssim
\sqrt{\mathcal{E}_\varepsilon (\nu)}.
$$
This implies that the set on which $n^{\varepsilon/4} \hat{\alpha}(\eta, n)$ is unbounded
must be of zero $\varepsilon$-capacity  (i.e., it does not support a measure of finite
$\varepsilon$-energy).

As the Hausdorff dimension of sets of zero $\varepsilon$-capacity is less than or equal
to $\varepsilon$ (see \cite[\S IV.1]{Carleson}), this completes the proof of the fact
that $S$ has zero Hausdorff dimension.

The second assertion follows since $d\mu_{{\rm sing}}$ is supported on the set $S$;
compare \cite[Corollary~10.8.4]{simon}.
\end{proof}

%%%%%%%%%%%%%%%%%%%%%%%%%%%%%%%%%%%%%%%%%%%%%%%%%%%%%%%%%%%%%%%%%%%%%%%%%%%%%%%%%%%%%
\section{Absence of a Singular Continuous Component}
%%%%%%%%%%%%%%%%%%%%%%%%%%%%%%%%%%%%%%%%%%%%%%%%%%%%%%%%%%%%%%%%%%%%%%%%%%%%%%%%%%%%%

In this section we employ ideas of Kiselev \cite{k} to show that there is no singular
continuous component when \eqref{vcass} holds. The preparatory work from the previous
section will be crucial.

The first step is to study the number of resonant points on the unit circle at which the
Pr\"ufer radius may be large. Using \eqref{vcass} and an almost-orthogonality lemma from
\cite{kls}, we will show that their number must be bounded by an explicit constant.

We first recall \cite[Lemma~4.4]{kls}:

\begin{lemma}\label{klslemma}
Let $e_1, \ldots, e_K$ be unit vectors in a Hilbert space $\mathcal{H}$ with
$$
Q = K \sup_{k \not= l} | \langle e_k, e_l \rangle | < 1.
$$
Then, for any $g \in \mathcal{H}$,
$$
\sum_{l = 1}^K | \langle g, e_l \rangle |^2 \le (1 + Q) \|g\|^2.
$$
\end{lemma}

Below, the Hilbert spaces in question will be given by $\mathcal{H}_n = \C^n$ with inner
product
$$
\langle f,g \rangle_{\mathcal{H}_n} = \sum_{j = 0}^{n-1} \overline{f(j)} g(j) (1+j).
$$

Recall Abel's formula (summation by parts), which reads
\begin{equation}\label{abelform}
\sum_{j = m}^n (\delta^+ a)(j) \cdot b(j) = a(n+1) \cdot b(n) - a(m) \cdot b(m-1) -
\sum_{j = m}^n a(j) \cdot (\delta^- b)(j).
\end{equation}
Here, $a, b$ are sequences, $(\delta^+ a)(j) = a(j+1) - a(j)$, and $(\delta^- b)(j) =
b(j) - b(j-1)$.

\begin{lemma}\label{abelstep}
Assume \eqref{vcass}. If $g$ is a real-valued sequence with $|(\delta^- g)(j)| \le B
|\alpha_{j-2}|$ for a suitable $B > 0$, then there are constants $C_1, C_2 > 0$ such
that, for $\xi \in (0,1)$, we have
$$
\sup_{n \ge 1} \left| \sum_{j = 1}^n j^{-1} e^{i[j\xi + g(j)]} \right| \le C_1 \log
(\xi^{-1}) + C_2.
$$
\end{lemma}

\begin{proof}
It suffices to consider $n > \pi \xi^{-1} + 1$. Let
$$
a(j) = - \sum_{k = j}^\infty k^{-1} e^{ik\xi} \quad \text{ and } b(j) = e^{ig(j)}.
$$
Note that
\begin{equation}\label{ajbound}
|a(j)| \lesssim \sum_{k = j}^{j + \lceil \pi \xi^{-1} \rceil} j^{-1} \lesssim \log \left(
1 + \frac{\pi \xi^{-1}}{j} \right).
\end{equation}
Clearly,
\begin{equation}\label{smallj}
\left| \sum_{j = 1}^{\lceil \pi \xi^{-1} \rceil} j^{-1} e^{\left( i[j \xi + g(j)]
\right)} \right| \lesssim \log( \pi \xi^{-1} ).
\end{equation}
On the other hand,
$$
\left| \sum_{j = \lceil \pi \xi^{-1} \rceil + 1}^{n} j^{-1} e^{\left( i[j \xi + g(j)]
\right)} \right| = \left| \sum_{j = \lceil \pi \xi^{-1} \rceil + 1}^{n} (\delta^+ a)(j)
\cdot b(j) \right|,
$$
which, by \eqref{abelform}, is equal to
$$
\left| a(n+1) \cdot b(n) - a(\lceil \pi \xi^{-1} \rceil + 1) \cdot b(\lceil \pi \xi^{-1}
\rceil) - \sum_{j = \lceil \pi \xi^{-1} \rceil + 1}^n a(j) \cdot (\delta^- b)(j) \right|.
$$
By \eqref{vcass}, \eqref{ajbound}, and the assumption on $g$, $|(\delta^- g)(j)| \le B
|\alpha_{j-2}|$, this expression is bounded by a constant only depending on $A$ and $B$.
(Split the sum into dyadic blocks, apply Cauchy-Schwarz, and then \eqref{vcass}.)
Combining this bound with \eqref{smallj}, the lemma follows.
\end{proof}

Write $A(n,\eta)$ for $A(n,\eta,0)$. We will consider situations where there are
$\eta_1,\ldots,\eta_K$ such that
\begin{equation}\label{larger}
| A(n,\eta_l) | \ge \frac{\log n}{14} \text{ for } l = 1,\ldots,K
\end{equation}
and
\begin{equation}\label{separation}
\min_{k \not= l} d(\eta_k,\eta_l) \ge n^{-1/(3K^2)}.
\end{equation}
Our goal is to bound $K$ from above. This is accomplished by the following lemma.

\begin{lemma}\label{l43}
Assume \eqref{vcass}. Then there are constants $n_0$ and $K_{{\rm max}} = K_{{\rm
max}}(A)$ such that for $n \ge n_0$, there can be no more than $K_{{\rm max}}$ points in
$[0,2\pi)$ for which \eqref{larger} and \eqref{separation} hold.
\end{lemma}

\begin{proof}
Consider $\eta_1,\ldots,\eta_K$ for which \eqref{larger} and \eqref{separation} hold.
Define the following vectors in $\mathcal{H}_n$:
$$
e_l(j) = E_n^{-1/2} e^{i[(j+1)\eta_l + 2 \theta_j(\eta_l,0)]} (1+j)^{-1}, \quad 1 \le l
\le K.
$$
The normalization constant $E_n$ is chosen so that the vectors have norm one. Obviously,
$E_n = \log n + O(1)$.

Now, for $k \not= l$, we have for $n$ large enough,
\begin{align*}
| \langle e_k, e_l \rangle | & = \frac{1}{E_n} \left| \sum_{j = 0}^{n-1} (1+j)^{-1}
e^{\left( - i[(j+1)\eta_k - 2 \theta_j(\eta_k,0)] + i[(j+1)\eta_l + 2 \theta_j(\eta_l,0)]
\right)} \right| \\
& = \frac{1}{E_n} \left| \sum_{j = 1}^{n} j^{-1} e^{\left( i[j (\eta_l -
\eta_k) - 2 \theta_{j-1}(\eta_k,0) + 2 \theta_{j-1}(\eta_l,0)] \right)} \right| \\
& \le \frac{C}{K^2}.
\end{align*}
We applied Lemma~\ref{abelstep} together with \eqref{separation} in the last step. The
constant $C$ depends only on $A$.

Thus, when $K > C$, we may apply Lemma~\ref{klslemma} and obtain for any $n \ge n_0$ and
$g \in \mathcal{H}_n$,
\begin{equation}\label{cons}
\sum_{l = 1}^K \left| \langle g, e_l \rangle_{\mathcal{H}_n} \right|^2 \le 2
\|g\|_{\mathcal{H}_n}^2.
\end{equation}
Let us apply \eqref{cons} to $g = (\alpha_0,\ldots,\alpha_{n-1})$. Due to \eqref{vcass},
the right-hand side can be estimated as follows:
$$
2 \|g\|_{\mathcal{H}_n}^2 = 2 \sum_{j = 0}^{n-1} |\alpha_j|^2 (j+1) \le 2 A \log n.
$$
On the other hand, by \eqref{larger},
$$
\left| \langle g, e_l \rangle_{\mathcal{H}_n} \right| = E_n^{-1/2} | A(n,\eta_l) | \ge
E_n^{-1/2} \frac{\log n}{14}.
$$
Consequently, \eqref{cons} implies that if $K > C$ and $n \ge n_0$,
$$
\frac{K ( \log n )^2}{196 E_n}  \le 2 A \log n.
$$
This shows that $K \le \tilde C$, with $\tilde C$ roughly being equal to $392 A$.
Therefore, we must have $K \le \max\{C,\tilde C \}$ whenever \eqref{larger} and
\eqref{separation} hold for $n \ge n_0$.
\end{proof}

Let us turn to the proof of the main theorem. Given the results above, we may from now on
follow the arguments of Kiselev in \cite{k} quite closely:

\begin{proof}[Proof of Theorem~\ref{main}.]
Assume that the singular continuous part of $d\mu$ is non-trivial. Fix an interval $I
\subset [0,2\pi)$ such that $\mu_{{\rm sc}}(I) = \Delta > 0$. Since $d \mu_{{\rm sc}}$ is
continuous, we can achieve that $\mu_{{\rm sc}}(J)$ is as small as we want if $J$ is any
subinterval of $I$ of sufficiently small length.

In particular, we can find $\varepsilon_0 \in (0,1)$ that satisfies the following
conditions ($K_{{\rm max}}$ and $n_0$ are the constants from Lemma~\ref{l43}):
\begin{itemize}

\item[(i)] $\lceil \varepsilon_0^{-3} \rceil > n_0$.

\item[(ii)] $\mu_{{\rm sc}}(J) < \frac{\Delta}{32 K_{{\rm max}}^3}$ for all intervals $J
\subseteq I$ with $|J| \le \varepsilon_0^{K_{{\rm max}}^{-2}}$.

\item[(iii)] $\frac{\varepsilon_0^{1/2}}{1 - \varepsilon_0^{1/2}} \le \frac{\Delta}{32
K_{{\rm max}}^3}$.

\item[(iv)] The last inequality holds in \eqref{muest} below.

\item[(v)] It is small enough so that we may obtain \eqref{larger} below.

\end{itemize}

We say that an interval $J \subset I$ belongs to \textit{scale} $m$ if $|J| =
\varepsilon_m := \varepsilon_0^m$. Two intervals of scale $m$ are called
\textit{separated} if the distance between their centers exceeds $3
\varepsilon_m^{K_{{\rm max}}^{-2}}$. An interval $J$ of scale $m$ is called
\textit{singular} if $\mu_{{\rm sc}}(J) > \varepsilon_m^{1/2}$.

We first show that there are no more than $K_{{\rm max}}$ separated singular intervals at
each scale. Assume that there are $K > K_{{\rm max}}$ separated singular intervals of
scale $m$: $J_1,\ldots,J_K$. Let $n_m = \lceil \varepsilon_m^{-3} \rceil$. Recall that
$d\mu_{n_m}$ denotes the Bernstein-Szeg\"o approximation of $d\mu$ at level $n_m$. Using
Lemma~\ref{bsaest}, we see that
\begin{equation}\label{muest}
\mu_{n_m}(3J_l) \ge \mu(J_l) - C \varepsilon_m > \varepsilon_m^{1/2} - C \varepsilon_m
\ge \tfrac12 \varepsilon_m^{1/2}.
\end{equation}
By \eqref{vcass} and \eqref{bsa},
$$
\frac{d\mu_{n_m}}{d\eta} (\eta) \sim R_{n_m}^{-2}(\eta,0).
$$
Thus, there are $\eta_l \in 3J_l$, $1 \le l \le K$, such that $R_{n_m}^{-2}(\eta,0)
\gtrsim \varepsilon_m^{-1/2}\!$, with a uniform implicit constant. In other words,
\eqref{larger} holds if $\varepsilon_0$ is small enough. Moreover, $\min_{k \not= l}
d(\eta_k,\eta_l) \ge \varepsilon_m^{K_{{\rm max}}^{-2}}$ because the intervals
$J_1,\ldots,J_K$ are separated. Thus, this yields a contradiction to Lemma~\ref{l43}.

Now write $S_m$ for the union of all singular intervals at scale $m$. This set can be
covered by at most $8K_{{\rm max}}$ intervals of size $\varepsilon_m^{K_{{\rm
max}}^{-2}}$, or else we can find more than $K_{{\rm max}}$ separated singular intervals
at scale $m$. By property~(ii) of $\varepsilon_0$, we get
$$
\mu_{{\rm sc}}(S_m) \le 8K_{{\rm max}} \times \frac{\Delta}{32 K_{{\rm max}}^3} =
\frac{\Delta}{4K_{{\rm max}}^2}
$$
for every $m \ge 1$. Now consider $m \ge K_{{\rm max}}^2$ and let $\tilde m = \lfloor
mK_{{\rm max}}^{-2} \rfloor \ge 1$. If $J_l^{(m)}$ is a singular interval at scale $m$
that obeys $\mu_{{\rm sc}}(J_l^{(m)}) > \varepsilon_{\tilde m}^{1/2}$, it must a subset
of $S_{\tilde m}$ since it can clearly be extended to a singular interval at scale
$\tilde m$. Thus, the set
$$
S_m \setminus \bigcup_{l < m} S_l
$$
can be covered by at most $8K_{{\rm max}}$ intervals of length $\varepsilon_m^{K_{{\rm
max}}^{-2}}$ and each of these intervals obeys $\mu_{{\rm sc}}(\cdot) \le
\varepsilon_{\tilde m}^{1/2}$. Consequently,
$$
\mu_{{\rm sc}}\left( S_m \setminus \bigcup_{l < m} S_l \right) \le 8K_{{\rm max}}
\varepsilon_{\tilde m}^{1/2}.
$$
Each $\tilde m$ corresponds to $K_{{\rm max}}^2$ values of $m$. Thus,
\begin{align*}
\mu_{{\rm sc}}\left( \bigcup_{m=1}^\infty S_m \right) & \le K_{{\rm max}}^2 \times
\frac{\Delta}{4K_{{\rm max}}^2} + K_{{\rm max}}^2 \times \sum_{\tilde m = 1}^\infty
8K_{{\rm max}} \varepsilon_{\tilde
m}^{1/2} \\
& = \frac{\Delta}{4} + \sum_{\tilde m = 1}^\infty 8K_{{\rm max}}^3 \varepsilon_0^{\tilde m /2} \\
& = \frac{\Delta}{4} +  8 K_{{\rm max}}^3 \frac{\varepsilon_0^{1/2}}{1 - \varepsilon_0^{1/2}} \\
& \le \frac{\Delta}{2}.
\end{align*}
In the last step, we used property (iii) of $\varepsilon_0$.

By zero-dimensionality (cf.~Proposition~\ref{p31}), $\mu_{{\rm sc}}|_{I}$ is supported by
the set
$$
D = \{ \eta \in I : \limsup_{ \delta \to 0 } \frac{\mu(k -
\delta,k+\delta)}{(2\delta)^{1/2}} = \infty \}.
$$
See, for example, \cite{rogers}. Thus, for each $k \in D$, there is a sequence $\delta_n
\to 0$ such that
$$
\frac{\mu(k - \delta_n,k+\delta_n)}{(2\delta_n)^{1/2}} \to \infty.
$$
For $n$ large, define $m_n$ by
$$
\frac{\varepsilon_{m_n}}{2} \ge \delta_n > \frac{\varepsilon_{m_n + 1}}{2} =
\frac{\varepsilon_0 \varepsilon_{m_n}}{2}.
$$
We obtain
$$
\frac{\mu(k - \varepsilon_{m_n}/2,k + \varepsilon_{m_n}/2)}{(\varepsilon_0
\varepsilon_{m_n})^{1/2}} \ge \frac{\mu(k - \delta_n,k+\delta_n)}{(2\delta_n)^{1/2}} \to
\infty.
$$
It follows that $k \in \bigcup_{m=1}^\infty S_m$ and hence
$$
0 < \Delta = \mu_{{\rm sc}}(I) = \mu_{{\rm sc}}(D) \le \mu_{{\rm sc}}\left(
\bigcup_{m=1}^\infty S_m \right) \le \frac{\Delta}{2},
$$
a contradiction.
\end{proof}

\end{document}